\documentclass[reqno,11pt]{amsart}
\usepackage{amssymb}
\usepackage{amsmath}
\usepackage{times}
\usepackage{graphics}
\usepackage{hyperref}

\newtheorem{theorem}{Theorem}[section]
\newtheorem{lemma}[theorem]{Lemma}
\newtheorem{proposition}[theorem]{Proposition}

\newtheorem{corollary}[theorem]{Corollary}

\begin{document}
\title{On eigenvalues of the sum of two random projections}
\author{V. Kargin }
\thanks{Statistical Laboratory, Department of Mathematics, University of
Cambridge, Cambridge, UK. e-mail: v.kargin@statslab.cam.ac.uk; vladislav.kargin@gmail.com}
\date{May 2012}
\maketitle

\begin{center}
\textbf{Abstract}
\end{center}

\begin{quotation}
We study the behavior of eigenvalues of matrix $P_N+Q_N$ where $P_N$ and $Q_N
$ are two $N$-by-$N$ random orthogonal projections. We relate the joint
eigenvalue distribution of this matrix to the Jacobi matrix ensemble and
establish the universal behavior of eigenvalues for large $N.$ The limiting
local behavior of eigenvalues is governed by the sine kernel in the bulk and
by either the Bessel or the Airy kernel at the edge depending on parameters.
We also study an exceptional case when the local behavior of eigenvalues of $%
P_N+Q_N$ is not universal in the usual sense.
\end{quotation}

\section{Introduction and Notation.}

Let $A$ and $B$ be two complex Hermitian matrices, perhaps random, and
consider matrix $X=A+UBU^{\ast }$ where $U$ is a uniformly-distributed
unitary matrix. What can be said about the distribution of eigenvalues of $X$
if the size of the matrices is large? Especially, what can be said about the
properties of this distribution on the local scale?

It is known that if $B$ is a random matrix with Gaussian entries, then the
eigenvalues of $A+B$ can be described as a determinantal process
(Brezin-Hikami \cite{brezin_hikami97}). Moreover, if in addition $A$ is a
Wigner matrix (that is, a Hermitian matrix with independent upper-diagonal
entries) and if the entries of $A$ and $B$ are comparable in size, then the
local properties of matrix $A+UBU^{\ast }$ are the same as if this matrix
were from GUE, Gaussian Unitary Ensemble, (Johansson \cite{johansson01}). On
the other hand if $B$ is a multiple of the identity matrix then the
eigenvalues of $A+UBU^{\ast }$ are obviously simply the shifted eigenvalues
of the matrix $A.$

The main question at stake here is that of the universality of the
eigenvalue distribution of $A+UBU^{\ast}$. The universality hypothesis
claims that for a large class of matrices the local correlations of
eigenvalues should approach certain specific functions when the size of
matrices grows. The choice of these functions should depend on the overall
symmetry properties of matrices such as whether their are complex Hermitian
or real symmetric, and on the properties of the global limit density of
eigenvalues. Roughly speaking, if the global density is positive and smooth
at a point $x$, then the limit local behavior at $x$ should be described by
the determinantal point process with the sine kernel. If the global density
has a square root singularity at $x$, then the local behavior around this
point should be described by the determinantal point process with the Airy
kernel, and so on.

The universality hypothesis originated in pioneering works of Dyson and
Mehta. (See, for example, the discussion in \cite{mehta04}.) It was verified
in various particular cases in \cite{soshnikov99}, \cite{soshnikov02}, \cite%
{bleher_its99}, \cite{dkmvz99}. Very recently, universality was proved for
the important cases of Wigner matrices and beta ensembles. See \cite%
{erdos_schlein_yau11}, \cite{tao_vu11} and a review in \cite{erdos_yau11}.

Now suppose that $A$ and $B$ are different from a multiple of the identity
matrix. Is it true that the behavior of eigenvalues of the matrix $%
A+UBU^{\ast }$ is universal?

The goal of this paper is to investigate this question in the case when $A$
and $B$ are random orthogonal projections, which we will denote $P$ and $Q$.
That is, $P$ and $Q$ are Hermitian (or real symmetric) matrices such that $%
P^{2}=P$ and $Q^{2}=Q.$ We will assume that the matrices $P$ and $Q$ are $N$%
-by-$N,$ and their ranks are $p$ and $q$ respectively. We will always
consider orthogonal projections in this paper so we omit the adjective
orthogonal.

We will call a projection matrix $P$ a \emph{random Hermitian (or real
symmetric) projection} if its entries are complex (or real) random variables
and their joint distribution is invariant with respect to unitary (or real
orthogonal) conjugations.

In a summary, our main finding is that in most situations the behavior of
eigenvalues of $P+Q$ for large $N$ is universal in the sense that it is
similar to the behavior of the eigenvalues of classical matrix ensembles.
The exceptional case is the local behavior of eigenvalues of $P+Q$ near $x=1$
when the ranks of $P$ and $Q$ are both approximately $N/2,$ and we discuss
this case below.

The reason why the eigenvalue distribution of $P+Q$ is easier to analyze
than that of the general case $A+B$ is that there is a relation of this
distribution to the distribution of the Jacobi ensemble of random points. We
formulate this relation in a slightly greater generality.

\begin{theorem}
\label{theorem_eigenvalue_distribution}Suppose $P$ is a random Hermitian (or
real symmetric) projection of rank $p,$ and $Q$ is an independent Hermitian
(or real symmetric) random projection of rank $q,$ $p\leq q.$ Assume $%
p+q\leq N$ and $\theta \neq 0.$ Then with probability $1,$ $P+\theta Q$ has $%
a\equiv q-p$ eigenvalues $\theta $ and $b\equiv N-p-q$ eigenvalues $0.$ The
density of the absolutely continuous part of the eigenvalue distribution of $%
P+\theta Q$ is induced by a transformation of the following densities: 
\newline
(i) 
\begin{equation}
f^{\left( 2\right) }\left( x_{1},\ldots ,x_{p}\right) =c_{2}\prod_{1\leq
i<j\leq p}\left\vert x_{i}-x_{j}\right\vert
^{2}\prod_{i=1}^{p}x_{i}^{a}\left( 1-x_{i}\right) ^{b},\text{ }x_{i}\in %
\left[ 0,1\right] ,  \label{ensemble_complex_Jacobi}
\end{equation}%
in the Hermitian case, and \newline
(ii) 
\begin{equation}
f^{\left( 1\right) }\left( x_{1},\ldots ,x_{p}\right) =c_{1}\prod_{1\leq
i<j\leq p}\left\vert x_{i}-x_{j}\right\vert \prod_{i=1}^{q}x_{i}^{a^{\prime
}}\left( 1-x_{i}\right) ^{b^{\prime }},\text{ }x_{i}\in \left[ 0,1\right] .
\label{ensemble_real_Jacobi}
\end{equation}%
in the real symmetric case, with $a^{\prime }=\left( a-1\right) /2$ and $%
b^{\prime }=\left( b-1\right) /2.$\newline
The transformation is two-fold and given by the formula 
\begin{equation*}
\lambda _{i}=\frac{1}{2}\left( 1+\theta \pm \sqrt{\left( 1-\theta \right)
^{2}+4\theta t_{i}}\right) .
\end{equation*}
\end{theorem}

This theorem allows us to reduce the study of spectral properties of $P+Q$
to the study of the corresponding Jacobi ensembles. In particular, consider
a sequence of matrices $P_{N},$ $Q_{N}$ with increasing size $N.$ We use
notation $p_{N},q_{N}$ to denote their ranks.

\begin{corollary}
Let $\mathcal{N}\left( I\right) $ denote the number of eigenvalues of the
Hermitian matrix $P_{N}+Q_{N}$ in an interval $I.$ Assume that $%
p_{N}/N\rightarrow p>0$ and $q_{N}/N\rightarrow q\geq p\,$as $N\rightarrow
\infty ,$ and suppose that $p+q\leq 1.$ If $I\subset (0,1)$ or $I\subset
(1,2)$ and $I$ has positive length$,$ then $\mathrm{Var}\left( \mathcal{N}%
\left( I\right) \right) \sim \pi ^{-2}\log N$ for large $N.$ Moreover, the
normalized random variable 
\begin{equation*}
\frac{\mathcal{N}\left( I\right) -\mathbb{E}\mathcal{N}\left( I\right) }{%
\sqrt{\mathrm{Var}\left( \mathcal{N}\left( I\right) \right) }}
\end{equation*}%
converges to a standard Gaussian random variable.
\end{corollary}

Note: The symbol $\subset $ denotes proper inclusion. A similar result holds
for real symmetric matrices $P_{N}+Q_{N}$ with a different value of the
limit variance.

\textbf{Proof of Corollary:} Let for definiteness $J=\left( \alpha ,\beta
\right) \subset \left( 0,1\right) $ and $I=(1+\alpha ^{2},1+\beta
^{2})\subset \left( 1,2\right) .$ Then $\mathcal{N}\left( I\right) =\mathcal{%
N}^{\prime }\left( J\right) $, where $\mathcal{N}^{\prime }\left( J\right) $
denote the number of points from ensemble (\ref{ensemble_complex_Jacobi}) in
interval $J.$ This ensemble is determinantal, hence we can apply the theorem
due to Costin, Lebowitz, and Soshnikov (\cite{costin_lebowitz95} and \cite%
{soshnikov00a} ) and conclude that 
\begin{equation*}
\frac{\mathcal{N}^{\prime }\left( J\right) -\mathbb{E}\mathcal{N}^{\prime
}\left( J\right) }{\sqrt{\mathrm{Var}\left( \mathcal{N}^{\prime }\left(
J\right) \right) }}
\end{equation*}%
converges to a standard Gaussian random variable, provided that $\mathrm{Var}%
\left( \mathcal{N}^{\prime }\left( J\right) \right) $ grows to infinity.

Under the assumption that $p_{N}/N\rightarrow p>0$ and $q_{N}/N\rightarrow
q\geq p\,$as $N\rightarrow \infty ,$ and that $p+q\leq 1,$ it is known that
the limit kernel for the corresponding Jacobi ensemble (\ref%
{ensemble_complex_Jacobi}) is the sine kernel. See, for example, \cite%
{collins05} or \cite{bbds06}. Hence the calculation of the variance can be
done similar to the calculation of variance in the case of the Gaussian
Unitary Ensemble (see calculations in \cite{gustavsson05} and related
formula (16.1.3) and Appendix A.38 in \cite{mehta04}). This calculation
results in the variance asymptotically equal to $\pi ^{-2}\log N$ for large $%
N.$ $\square $

Similarly, the other properties of eigenvalues of the Hermitian matrix $%
P_{N}+Q_{N}$ in the interval $J$ are analogous to the properties of the
points from the Jacobi ensemble for the simple reason that the
transformation $\lambda \rightarrow \left( \lambda -1\right) ^{2}$ is
locally linear on this interval. In particular, in the limit these
properties can be described by the sine kernel.

For the behavior at the edge of the support, the results are also similar to
the case of the classical ensembles and this behavior can be described by
either the Bessel kernel (if $p+q=1$) or the Airy kernel (if $p+q<1$). The
theorem below gives details for the case $p+q<1$.

\begin{theorem}
\label{thm_edge_behavior} Let $P_{N}$ and $Q_{N}$ be two random projections
in $\mathbb{C}^{N}$ with ranks $p_{N}$ and $q_{N}$ respectively. Assume that 
$\lim_{N\rightarrow \infty }\frac{p_{N}}{N}=p$ and $\lim_{N\rightarrow
\infty }\frac{q_{N}}{N}=q$, and suppose that $p+q<1.$ \ Let $\lambda
_{\left( 1\right) }$ denotes the largest eigenvalue of $P_{N}+Q_{N}$ which
is not equal to $2$ and let 
\begin{equation*}
\mu =\sqrt{q(1-p)}+\sqrt{p\left( 1-q\right) }.
\end{equation*}%
Then there exists $\sigma >0$ such that for every $u$ 
\begin{equation*}
\mathbb{P}\left\{ \frac{\lambda _{\left( 1\right) }-\mu }{\sigma N^{-2/3}}%
\leq u\right\} \rightarrow F_{2}\left( u\right) ,
\end{equation*}%
where $F_{2}\left( u\right) $ is the Tracy-Widom distribution function for
the edge of the Gaussian unitary ensemble.
\end{theorem}

This theorem directly follows from recent results about the distribution of
the largest eigenvalue in the Jacobi ensemble with changing parameters. A
similar result holds for real symmetric projections except that $\sigma $ is
different and the probability converges to $F_{1}\left( u\right) $ instead
of $F_{2}\left( u\right) .$

The only unusual situation occurs at $x=1$ when this point belongs to the
support of the absolutely continuous part of the limit eigenvalue
distribution, which occurs when $p=q=1/2.$ We will analyze this special
situation in the following theorem.

Let $J_{a}\left( z\right) $ denote the Bessel function with index $a.$

\begin{theorem}
\label{theorem_near_lambda_1}Let $P_{N}$ and $Q_{N}$ be two random
projections in $\mathbb{C}^{N}$ with ranks $p_{N}$ and $q_{N}$ respectively.
Assume that $q_{N}-p_{N}=a\geq 0,$ and $N-p_{N}-q_{N}=b\geq 0.$ Let $%
s_{N}=t/\left( \sqrt{2}p_{N}\right) ,$ where $t>0.$ Then the expected number
of eigenvalues of $P+Q$ in the interval $I_{N}=[1,1+s_{N}]$ is 
\begin{equation*}
\mathbb{E}\mathcal{N}\left( I_{N}\right) =\int_{0}^{t}x\left( J_{a}\left(
x\right) ^{2}-J_{a+1}\left( x\right) J_{a-1}\left( x\right) \right) dx+o(1)
\end{equation*}%
for large $N.$
\end{theorem}

By using the power series for the Bessel functions, 
\begin{equation*}
J_{a}\left( z\right) =\sum_{m=0}^{\infty }\frac{\left( -1\right) ^{m}\left(
z/2\right) ^{2m+a}}{m!(m+a)!},
\end{equation*}%
we  obtain the following corollary.

\begin{corollary}
\label{corollary_number_eigenvalues}For large $N,$%
\begin{equation*}
\mathbb{E}\mathcal{N}\left( I_{N}\right) \sim \frac{2}{\left( a+1\right)
!^{2}}\left( \frac{t}{2}\right) ^{2a+2}+R\left( t\right) ,
\end{equation*}%
where $R\left( t\right) $ is a differentiable function such that $%
\lim_{t\rightarrow 0}R(t)/t^{2a+4}$ exists and finite.
\end{corollary}

(Here the asymptotic equivalence  $\sim $ is understood relative to increase
in $N.$) For example, if $a=0,$ then $\mathbb{E}\mathcal{N}\left(
I_{N}\right) \sim t^{2}/2+R(t)$ where $\lim_{t\rightarrow 0}R\left( t\right)
/t^{4}<\infty .$

In order to appreciate the significance of this corollary, note that under
assumptions of Theorem \ref{theorem_near_lambda_1}, the eigenvalue
distribution of $P_{N}+Q_{N}$ weakly converges (in probability) to a fixed
probability distribution which is absolutely continuous on $\left(
0,2\right) $ and has the density 
\begin{equation*}
f\left( \lambda \right) =\frac{1}{\pi }\frac{1}{\sqrt{\lambda (2-\lambda )}}.
\end{equation*}%
This can be shown either by using a relation of this ensemble with free
probability theory (see \cite{pastur_vasilchuk00}) or by computing the limit
distribution for the corresponding Jacobi ensemble.

In particular, the limiting density is not zero at $\lambda =0$ and does not
depend on $a$ or $b.$ Hence we can see that the expected number of
eigenvalues in the interval $I_{N}$ cannot be estimated as $N\times $ length
of $I_{N}\times $ density + terms of smaller order, which would result in a
term of the order $ct$ with $c>0.$ For the ensemble $P_{N}+Q_{N},$ the
expected number of eigenvalues in $I_{N}$ is about $t^{2}/2$ and
significantly smaller than this estimate.

Thus, the local behavior near $x=1$ differs strongly from the local behavior
of the GUE eigenvalues. In the case of GUE, a similarly scaled limit gives $%
\mathbb{E}\mathcal{N}\left( I_{N}\right) \sim \rho \left( x\right) t$ where $%
\rho \left( x\right) $ is the density of the semicircle law. Hence, the
result in Corollary \ref{corollary_number_eigenvalues} is in contradiction
with the heuristic that the local behavior of eigenvalues of an Hermitian
ensemble is similar to that of GUE eigenvalues provided that the limit
eigenvalue density is non-zero and smooth at a given point. (For example, at
the beginning of Section 6.3.3 on p. 113 of \cite{kuijlaars11}, it is stated
that \textquotedblleft In the context of Hermitian matrix models, the
universality may be stated as the fact that the global eigenvalue regime
determines the local eigenvalue regime.\textquotedblright )

A possible intuitive explanation of this phenomenon is the additional
symmetry of the model. This symmetry forces eigenvalues either land exactly
on $x=1,$ or come in pairs located symmetrically around $x=1.$ In both cases
the repulsion between eigenvalues ensures that for small $t$ the interval $%
(0,t/2p_{N})$ has fewer eigenvalues than the density function predicts.

Another interesting feature of the result in Theorem \ref%
{theorem_near_lambda_1} is that the limit does not depend on $b,$ that is,
on the excess of the dimension $N$ of the ambient space over the total rank
of $P$ and $Q.$ The only provision is that $b$ stays constant as $N$ grows.

Hence the local behavior of eigenvalues near $x=1$ exhibits universality
with respect to parameter $b$ but not with respect to parameter $a.$ The
difference between ranks of $P$ and $Q$ does matter for the limiting local
behavior of eigenvalues near $x=1$ despite the fact that it does not
influence the global limiting distribution.

For the case of real symmetric matrices we have an analogous result.

\begin{theorem}
\label{theorem_near_lambda_1_real}Let $P_{N}$ and $Q_{N}$ be two random
projections in $\mathbb{R}^{N}$ with ranks $p_{N}$ and $q_{N}$ respectively.
Assume that $q_{N}-p_{N}=a\geq 0,$ and $N-p_{N}-q_{N}=b\geq 0.$ Let $%
s_{N}=t/\left( \sqrt{2}p_{N}\right) ,$ where $t>0.$ Then for large $N$ the
expected number of eigenvalues of $P+Q$ in the interval $I_{N}=[1,1+s_{N}]$
is 
\begin{eqnarray*}
\mathbb{E}\mathcal{N}\left( I_{N}\right) &=&\int_{0}^{t}x\left( J_{a^{\prime
}+1}\left( x\right) ^{2}-J_{a^{\prime }+2}\left( x\right) J_{a^{\prime
}}\left( x\right) \right) dx \\
&&+\int_{0}^{t}J_{a^{\prime }+1}\left( x\right) \left(
\int_{0}^{x}J_{a^{\prime }+1}\left( u\right) du-1\right) dx+o(1),
\end{eqnarray*}%
where $a^{\prime }=\left( a-1\right) /2.$
\end{theorem}

By using the power series expansion for the Bessel functions we obtain the
following consequence.

\begin{corollary}
For large $N,$%
\begin{equation*}
\mathbb{E}\mathcal{N}\left( I_{N}\right) \sim \frac{1}{\left( a+3\right) 2^{%
\frac{a-1}{2}\Gamma \left( \frac{a+1}{2}\right) }}t^{\frac{a+3}{2}}+R\left(
t\right) .
\end{equation*}%
where $R\left( t\right) $ is a differentiable function such that $%
\lim_{t\rightarrow 0}R(t)/t^{\left( a+7\right) /2}$ exists and finite.
\end{corollary}

For example for $a=0,$ $\mathbb{E}\mathcal{N}\left( I_{N}\right) \sim
1/\left( 3\sqrt{\pi /2}\right) t^{3/2}+R(t),$ where $\lim_{t\rightarrow
0}R(t)/t^{7/2}<\infty .$

The rest of the paper organized as follows: we prove Theorem \ref%
{theorem_eigenvalue_distribution} in Section \ref%
{section_proof_of_main_theorem} after proving some results about the
relation of eigenvalues of $P+Q$ and of $PQP$ in Section \ref%
{section_eigenvalues_and_sing_values}. Section \ref{section_edge_behavior}
is about the (universal) edge behavior of eigenvalues. Section \ref%
{sections_eigenvalues_near_1} proves Theorems \ref{theorem_near_lambda_1}
and \ref{theorem_near_lambda_1_real} about the exceptional case when the
behavior of eigenvalues is not universal. And Section \ref%
{section_conclusion} concludes.

\section{Relation of eigenvalues of P+Q and PQP}

\label{section_eigenvalues_and_sing_values}

\begin{lemma}
Let $P$ and $Q$ be orthogonal projections and assume that $\mathrm{Range}%
P\cap \mathrm{Range}Q=0.$ Let $t_{i}>0$ be an eigenvalue of operator $PQP$.
Then, 
\begin{equation*}
\lambda _{i}=\frac{1}{2}\left( 1+\theta \pm \sqrt{\left( 1-\theta \right)
^{2}+4\theta t_{i}}\right)
\end{equation*}
are eigenvalues of $P+\theta Q.$
\end{lemma}

\textbf{Proof:} Let $u_{i}$ be an eigenvector of $PQP$ with eigenvalue $%
t_{i} $ and define $v_{i}:=QPu_{i}$. Then, we have%
\begin{equation*}
QPu_{i}=v_{i},\text{ and }Pv_{i}=t_{i}u_{i}.
\end{equation*}%
Clearly, $u_{i}\in \mathrm{Range}P,$ $v_{i}\in \mathrm{Range}Q.$ Hence, $%
Pu_{i}=u_{i}$ and $Qv_{i}=v_{i}.$

Define 
\begin{equation*}
\alpha _{i}=\frac{1}{2t_{i}}\left( -1+\theta \pm \sqrt{\left( 1-\theta
\right) ^{2}+4\theta t_{i}}\right) .
\end{equation*}%
Then,%
\begin{eqnarray*}
\left( P+\theta Q\right) \left( u_{i}+\alpha _{i}v_{i}\right)
&=&u_{i}+\alpha _{i}t_{i}u_{i}+\theta v_{i}+\theta \alpha _{i}v_{i}=\left(
1+\alpha _{i}t_{i}\right) \left( u_{i}+\frac{\theta \left( 1+\alpha
_{i}\right) }{\alpha _{i}\left( 1+\alpha _{i}t_{i}\right) }\alpha
_{i}v_{i}\right) \\
&=&\lambda _{i}\left( u_{i}+\alpha _{i}v_{i}\right) .
\end{eqnarray*}%
Note that $u_{i}+\alpha _{i}v_{i}\neq 0$ by the assumption that $\mathrm{%
Range}P\cap \mathrm{Range}Q=0.$ Hence $\lambda _{i}$ is an eigenvalue of $%
\left( P+\theta Q\right) .$ $\square $

\begin{lemma}
Suppose $P$ is a random projection of rank $p,$ and $Q$ is a projection of
rank $q,$ $p\leq q.$ Assume $p+q\leq N$ and $\theta \neq 0.$ Then, with
probability $1,$ the matrix $P+\theta Q$ has $2p$ distinct eigenvalues 
\begin{equation*}
\lambda _{i}^{\pm }=\frac{1}{2}\left( 1+\theta \pm \sqrt{\left( 1-\theta
\right) ^{2}+4\theta t_{i}}\right)
\end{equation*}%
where $t_{1},\ldots ,t_{p}$ are non-zero eigenvalues of operator $PQP,$ $%
t_{i}\in \left( 0,1\right) .$ In addition, with probability $1,$ $P$ has $%
q-p $ eigenvalues $\theta $ and $N-p-q$ eigenvalues $0.$
\end{lemma}

\textbf{Proof:} Let $[L]$ denote the dimension of linear space $L.$ Then, by
using the fact that with probability $1,$ the range and kernel spaces of $P$
and $Q$ are in general position relative to each other, we can calculate:%
\begin{eqnarray}
\left[ \mathrm{Ker}P\cap \mathrm{Ker}Q\right]  &=&N-p-q, \\
\lbrack \mathrm{Range}P\cap \mathrm{Range}Q] &=&\left( p+q-N\right) _{+}=0,
\label{range_intersection} \\
\lbrack \mathrm{Ker}P\cap \mathrm{Range}Q] &=&\left( N-p+q-N\right) _{+}=q-p,
\\
\lbrack \mathrm{Range}P\cap \mathrm{Ker}Q] &=&\left( p+N-q-N\right) _{+}=0,
\\
\lbrack \mathrm{Ker}QP] &=&[\mathrm{Range}P\cap \mathrm{Ker}Q]+[\mathrm{Ker}%
P] \\
&=&N-p, \\
\lbrack \mathrm{Range}QP] &=&N-[\mathrm{Ker}QP]=p,
\end{eqnarray}%
which implies that $QP$ has $p$ positive singular values, and therefore $PQP$
has $p$ positive eigenvalues. In addition, (\ref{range_intersection})
implies that all singular values of $QP$ are smaller than $1.$

Therefore, by previous lemma $Q+P$ has at least $2p$ eigenvalues (possibly
counting with multiplicities) that are distinct from $0$ and $\theta .$

With probability $1$ the singular values of $QP$ in the interval $\left(
0,1\right) $ are distinct. Indeed, the set $\mathcal{X}$ of $P$ such that $%
PQP$ has a multiple eigenvalue $t$ with $t\in \left( 0,1\right) $ is
characterised by an algebraic condition and therefore it is algebraically
closed. The set of all $P$ is algebraically irreducible. Hence $\mathcal{X}$
is either coincide with the set of all $P$ or has a smaller dimension than
the set of all $P$ and therefore has measure zero. It is possible to
construct an example such that all eigenvalues in the interval $\left(
0,1\right) $ are distinct. Therefore the measure of $\mathcal{X}$ is zero.

Therefore, the corresponding eigenvalues of $Q+P$ are distinct with
probability $1.$

In addition $\left[ \mathrm{Ker}P\cap \mathrm{Ker}Q\right] $ is an
eigenspace of $Q+P$ with dimension $N-p-q,$ and eigenvalue $0,$ and $[%
\mathrm{Ker}P\cap \mathrm{Range}Q]$ is an eigenspace of $Q+P$ with dimension 
$q-p$ and eigenvalue $\theta .$ By counting multiplicities we conclude that
we found all eigenvalues. $\square $

\section{Relation of eigenvalues of P+Q and the Jacobi ensemble}

\label{section_proof_of_main_theorem}

By using known results about the eigenvalues of $PQP$, we are able to prove
Theorem \ref{theorem_eigenvalue_distribution}. (These eigenvalues have been
studied in \cite{collins05} in connection with free probability.)

\textbf{Proof of Theorem \ref{theorem_eigenvalue_distribution}:} By a change
of coordinates we can assume without loss of generality that $Q$ is a
diagonal matrix that has $Q_{ii}=1$ for $i\leq q$ and $Q_{ii}=0$ for $i>q.$

We will further assume that $P$ is a projection on the span of $p$
independent random vectors $v_{1},\ldots ,v_{p}.$ The components of a vector 
$v_{i}$ are independent standard Gaussian variables (real in the symmetric
case and complex in the Hermitian case). If $X$ is an $N$-by-$p$ matrix with
columns $v_{1},\ldots ,v_{p},$then we can write: 
\begin{equation*}
X=\left[ 
\begin{array}{c}
A \\ 
B%
\end{array}%
\right] ,
\end{equation*}%
where $A$ is a $q$-by-$p$ matrix and $B$ is a $\left( N-q\right) $-by-$p$
matrix. After orthonormalization, we get 
\begin{equation}
Y=\left[ 
\begin{array}{c}
A \\ 
B%
\end{array}%
\right] \left( A^{^{\prime }}A+B^{\prime }B\right) ^{-1/2},
\label{formula_Y}
\end{equation}%
and the projection on the span of $v_{1},\ldots ,v_{p}$ is 
\begin{equation*}
P=YY^{\prime }.
\end{equation*}%
Formula (\ref{formula_Y}) implies that the non-zero $q$-by-$q$ block of $QPQ$
equals $A\left( A^{^{\prime }}A+B^{\prime }B\right) ^{-1}A^{\prime }.$
Hence, the non-zero eigenvalues of $PQP$ correspond to non-zero eigenvalues
of a $p$-by-$p$ matrix $M=\left( A^{^{\prime }}A+B^{\prime }B\right)
^{-1}A^{\prime }A$ \ It is known that the joint distribution of eigenvalues
of $M$ follows the law of the so-called Jacobi ensemble with parameters $%
\left( p,N-q,q\right) $. (See Section 3.6 and 7.2.5 of \cite{forrester10}
for more information about the Jacobi ensemble.) Namely, the density for the
complex case is 
\begin{equation*}
f^{\left( 2\right) }\left( x_{1},\ldots ,x_{p}\right) =c_{2}\prod_{1\leq
i<j\leq p}\left\vert x_{i}-x_{j}\right\vert
^{2}\prod_{i=1}^{p}x_{i}^{a_{2}}\left( 1-x_{i}\right) ^{b_{2}},\text{ }%
x_{i}\in \left[ 0,1\right] ,
\end{equation*}%
where $a_{2}=q-p$ and $b_{2}=N-p-q$ (see Section 8 in \cite{james64}).

The density for the real case is 
\begin{equation*}
f^{\left( 1\right) }\left( x_{1},\ldots ,x_{p}\right) =c_{1}\prod_{1\leq
i<j\leq p}\left\vert x_{i}-x_{j}\right\vert
\prod_{i=1}^{q}x_{i}^{a_{1}}\left( 1-x_{i}\right) ^{b_{1}},\text{ }x_{i}\in %
\left[ 0,1\right] .
\end{equation*}%
where $a_{1}=\left( q-p-1\right) /2$ and $b_{1}=\left( N-q-p-1\right) /2,$ .
(see Muirhead, Thm 3.3.4). This completes the proof. $\square $

\section{Eigenvalues at the edge}

\label{section_edge_behavior}

In the first step of this section we establish the limit of the eigenvalue
distribution on the global scale. In particular we identify the support of
the limit distribution. Then we study the behavior of the eigenvalues at the
edge by using the connection between eigenvalues of $P_{N}+Q_{N}$ and the
Jacobi ensemble.  We rely on  known results about the edge behavior of
eigenvalues from the Jacobi ensemble. As an aside, we derive the limit
global density for the Jacobi ensemble with variable parameters.

We abuse notation a bit and denote $\lim_{N\rightarrow \infty }\frac{p_{N}}{N%
}$ by $p$ and $\lim_{N\rightarrow \infty }\frac{q_{N}}{N}$ by $q$. That is,
in this section $p$ and $q$ denote the limits of fractional ranks of
matrices $P_{N}$ and $Q_{N}.$ We also use the notation 
\begin{eqnarray*}
a &:&=q-p\equiv \lim_{N\rightarrow \infty }\frac{a_{N}}{N}, \\
b &:&=1-p-q\equiv \lim_{N\rightarrow \infty }\frac{b_{N}}{N},\text{ and } \\
c &:&=p\left( 1-q\right) +q\left( 1-p\right) \equiv (1+a^{2}-b^{2})/2.
\end{eqnarray*}

Let $\lambda _{i}$ be ordered eigenvalues of $P_{N}+Q_{N},$ counted with
multiplicity. Let 
\begin{equation*}
\mu _{N}:=\frac{1}{N}\sum_{i=1}^{N}\delta _{\lambda _{i}};
\end{equation*}%
that is, $\mu _{N}$ is the eigenvalue distribution of $P_{N}+Q_{N}.$

\begin{proposition}
\label{prop_limit_of_pq}With probability $1$ the eigenvalue distribution of $%
P_{N}+Q_{N}$ weakly converges to a probability measure $\mu .$ The
absolutely continuous part of this measure is supported on $I_{1}\cup I_{2},$
where $I_{1}=1+I,$ $I_{2}=1-I,$ and 
\begin{equation*}
I=\left[ \left\vert \sqrt{q(1-p)}-\sqrt{p\left( 1-q\right) }\right\vert ,%
\sqrt{q(1-p)}+\sqrt{p\left( 1-q\right) }\right] .
\end{equation*}%
The density of this part is given by 
\begin{equation*}
\rho \left( x\right) =\frac{1}{\pi }\frac{\sqrt{-a^{2}+2c\left( x-1\right)
^{2}-\left( x-1\right) ^{4}}}{x(x-1)(x-2)}.
\end{equation*}%
In addition, $\mu $ has an atom with weight $\left\vert a\right\vert $ at $%
x=1.$ If $b>0$ then $\mu $ has an atom with weight $b$ at $x=0;$ and if $b<0$
then $\mu $ has an atom with weight $-b$ at $x=2.$
\end{proposition}

\textbf{Proof:} This result is a particular case of a theorem in \cite%
{pastur_vasilchuk00}, with the limit distribution computed by using the
standard techniques from free probability. $\square $

As a consequence, we can write down the limit level distribution for the
Jacobi ensemble with a changing weight function. (This limit level
distribution is known and goes back to the research of Wachter in \cite%
{wachter80}. However, we obtain it by using free probability calculations
and this seems to be knew.)

Let 
\begin{equation}
f^{\left( 2\right) }\left( x_{1},\ldots ,x_{N}\right) =c_{2}\prod_{1\leq
i<j\leq N}\left\vert x_{i}-x_{j}\right\vert
^{2}\prod_{i=1}^{N}x_{i}^{s_{N}}\left( 1-x_{i}\right) ^{t_{N}},\text{ }%
x_{i}\in \left[ 0,1\right] ,  \label{Hermitian_Jacobi_ensemble}
\end{equation}

Recall that the level density of this distribution is simply the marginal
distribution: 
\begin{equation*}
p_{N}^{\left( 2\right) }\left( x;s_{N},t_{N}\right) =\int f^{\left( 2\right)
}\left( x,x_{2},\ldots ,x_{N}\right) dx_{2}\ldots dx_{N}.
\end{equation*}

\begin{proposition}
\label{prop_limit_of_Jacobi}Let $p_{N}^{\left( 2\right) }\left(
x;s_{N},t_{N}\right) $ is the level density of the Hermitian Jacobi ensemble
with parameters $s_{N}$ and $t_{N}.$ Assume that $s_{N}/N\rightarrow s\geq 0$
and $t_{N}/N\rightarrow t\geq 0$ as $N\rightarrow \infty .$ Define: 
\begin{eqnarray*}
a &:&=\frac{s}{2+t+s}, \\
c &:&=\frac{1}{2}\left( 1+\frac{s^{2}-t^{2}}{\left( 2+t+s\right) ^{2}}%
\right) ,\text{ and} \\
d &:&=\frac{1}{\left( 2+t+s\right) ^{2}}\sqrt{\left( 1+s\right) \left(
1+t\right) \left( 1+t+s\right) }.
\end{eqnarray*}%
Then the level density $p_{N}^{\left( 2\right) }\left( x;s_{N},t_{N}\right) $
converges to the limit $p^{\left( 2\right) }\left( x;s,t\right) $ which is
supported on the interval 
\begin{equation*}
I=\left[ c-d,c+d\right] .
\end{equation*}%
Its density is%
\begin{equation*}
\rho ^{(2)}\left( x\right) =\frac{2+t+s}{2\pi }\frac{\sqrt{-a^{2}+2cx-x^{2}}%
}{x(x-1)}.
\end{equation*}
\end{proposition}

\textbf{Proof:} We use the correspondence between the Jacobi ensemble and
the eigenvalues of the sum of random projections $P_{N}+Q_{N}$ established
in Theorem \ref{theorem_eigenvalue_distribution}. For this we use projectios 
$P_{N}$ and $Q_{N}$ with ranks $p_{N}$ and $q_{N}$ and make sure that $%
p_{N}/N$ and $q_{N}/N$ approach $p$ and $q,$ which are defined as follows: 
\begin{equation*}
p:=\frac{1}{2+t+s}\text{ and }q:=\frac{1+s}{2+t+s}.
\end{equation*}%
By using Theorem \ref{theorem_eigenvalue_distribution}, it is easy to check
that the limit eigenvalue density of this matrix ensemble must coincide with
the limit level density of the Jacobi ensemble. Then, the conclusion of the
theorem follows from Proposition \ref{prop_limit_of_pq}. $\square $

Now let us find out what is the edge behavior of eigenvalues of $%
P_{N}+Q_{N}. $ The edge behavior of the Jacobi ensemble with variable
parameters was analysed in \cite{collins05} and \cite{johnstone08}. (The
case with the fixed parameters was studied in \cite{forrester93} and \cite%
{nagao_forrester95}.)

In particular, Johnstone found that the distribution of the largest
eigenvalue is concentrated near a point $x,$ which he defines as follows.

Let 
\begin{equation*}
\cos \phi =\frac{t-s}{2+t+s}\text{ and }\cos \gamma =\frac{t+s}{2+t+s}.
\end{equation*}%
Then, 
\begin{equation*}
x=\frac{1-\cos \left( \phi +\gamma \right) }{2}.
\end{equation*}%
(We change Johnstone's notations a bit to adjust them to our situation since
we consider the Jacobi ensemble on $\left[ 0,1\right] $ instead of $\left[
-1,1\right] .$ In addition, we omitted some terms of order $O\left(
N^{-1}\right) $ which are important for statistical applications but
irrelevant from the asymptotic point of view.)

An easy verification shows that $x$ coincides with the upper bound of the
limiting distribution support which we derived in Proposition \ref%
{prop_limit_of_Jacobi}, that is,  $x=c+d.$

Johnstone proved the following result (see Theorem 3 in \cite{johnstone08}
and a related result in Theorem 4.15 of \cite{collins05}):

\begin{theorem}[Johnstone]
\label{thm_johnstone}Assume that $t>0$ and let $x_{\left( 1\right) }$ denote
the largest point from $N$ points sampled according to distribution (\ref%
{Hermitian_Jacobi_ensemble}). Let  $s_{N}/N\rightarrow s\geq 0$ and $%
t_{N}/N\rightarrow t\geq 0$ as $N\rightarrow \infty .$ Then, there exists $%
\sigma >0,$ such that for every $u,$ 
\begin{equation*}
\mathbb{P}\left\{ \frac{\tanh ^{-1}\left( x_{\left( 1\right) }\right) -\tanh
^{-1}\left( x\right) }{\sigma N^{-2/3}}\leq u\right\} \rightarrow
F_{2}\left( u\right) ,
\end{equation*}%
where $F_{2}\left( u\right) $ is the Tracy-Widom distribution function for
the edge of the Gaussian unitary ensemble.
\end{theorem}

Johnstone gives an explicit expression for $\sigma $ and an explicit bound
for the error in this convergence result. He uses a slightly different
definition of $x$ which improves the second order convergence properties.
Also note that while the inverse hyperbolic tangent improves the convergence
properties for relatively small $N$, it is irrelevant from the asymptotic
point of view. A similar theorem was proved by Johnstone for the real
symmetric Jacobi ensemble, with a different $\sigma $ and $F_{1}\left(
u\right) $ instead of $F_{2}\left( u\right) .$

\textbf{Proof of Theorem \ref{thm_edge_behavior}:} This result follows
immediately from Theorems \ref{theorem_eigenvalue_distribution} and \ref%
{thm_johnstone}. $\square $

\section{Eigenvalues of P+Q near x=1}

\label{sections_eigenvalues_near_1}

\subsection{Hermitian case}

The density%
\begin{equation}
f\left( x_{1},\ldots ,x_{n}\right) =c\prod_{1\leq i<j\leq n}\left\vert
x_{i}-x_{j}\right\vert ^{2}\prod_{i=1}^{n}x_{i}^{a}\left( 1-x_{i}\right)
^{b},\text{ }x_{i}\in \left[ 0,1\right]
\end{equation}%
can be written as a determinant: 
\begin{eqnarray*}
f\left( x_{1},\ldots ,x_{n}\right) &=&\frac{1}{N!}\det \left( K_{n}\left(
x_{i},x_{j}\right) \right) _{i,j=1...n}; \\
K_{n}\left( x,y\right) &=&\sum_{k=0}^{n-1}Q_{k}\left( x\right) Q_{k}\left(
y\right) \\
&=&c_{n}\frac{Q_{n}\left( x\right) Q_{n-1}\left( y\right) -Q_{n-1}\left(
x\right) Q_{n}\left( y\right) }{x-y},
\end{eqnarray*}%
where $Q_{k}\left( x\right) $ are orthonormal polynomials with respect to
the weight $x^{a}\left( 1-x\right) ^{b}$ on the interval $\left[ 0,1\right] $
and $c_{n}$ is an appropriate constant.

\begin{lemma}
\label{lemma_scaling_limit}We have 
\begin{equation*}
\lim_{n\rightarrow \infty }\frac{1}{2n^{2}}K_{n}\left( \frac{t^{2}}{2n^{2}},%
\frac{t^{2}}{2n^{2}}\right) =\frac{1}{2}\left( J_{a}\left( t\right)
^{2}-J_{a+1}\left( t\right) J_{a-1}\left( t\right) \right) ,
\end{equation*}%
where $J_{a}\left( t\right) $ denote the Bessel function of index $a.$
\end{lemma}

\textbf{Proof:} This is a consequence of a remark in Section 7.2.5 and
Proposition 7.2.1 in \cite{forrester10}. $\square $

\textbf{Proof of Theorem \ref{theorem_near_lambda_1}:} By using the relation
between the eigenvalues of $P_{N}+Q_{N}$ and of the Jacobi ensemble, we have 
$\mathbb{E}\mathcal{N}\left( I_{N}\right) =\mathbb{E}\mathcal{N}^{\prime
}\left( J_{N}\right) ,$ where $\mathbb{E}\mathcal{N}^{\prime }\left(
J_{N}\right) $ is the expected number of eigenvalues of the Jacobi ensemble
in the interval $J_{N}=[0,t_{0}^{2}/(2p_{N}^{2})].$ Hence, by using Lemma %
\ref{lemma_scaling_limit}, we can write: 
\begin{eqnarray*}
\mathbb{E}\mathcal{N}\left( I_{N}\right)
&=&\int_{0}^{t_{0}^{2}/(2p_{N}^{2})}K_{p_{N}}\left( x,x\right) dx \\
&=&\int_{0}^{t_{0}}\frac{t}{p_{N}^{2}}K_{p_{N}}\left( \frac{t^{2}}{2p_{N}^{2}%
},\frac{t^{2}}{2p_{N}^{2}}\right) dt \\
&=&\int_{0}^{t_{0}}t\left( J_{a}\left( t\right) ^{2}-J_{a+1}\left( t\right)
J_{a-1}\left( t\right) \right) dt+o(1).
\end{eqnarray*}%
$\square $

\subsection{Real symmetric case}

In the real symmetric case the density of eigenvalue distribution and all
correlation functions (i.e., the marginals of this density) can be written
using quaternion determinants: 
\begin{equation*}
f\left( x_{1},\ldots ,x_{n}\right) =\frac{1}{N!}\mathrm{Q}\det \left(
K_{n}\left( x_{i},x_{j}\right) \right) _{i,j=1...n},
\end{equation*}%
where $K_{n}\left( x_{i},x_{j}\right) $ is a quaternion function. The matrix
representation of this quaternion is 
\begin{equation*}
K_{n}\left( x,y\right) =\left( 
\begin{array}{cc}
S_{n}\left( x,y\right)  & \widetilde{I}_{n}\left( x,y\right)  \\ 
D_{n}\left( x,y\right)  & S_{n}\left( y,x\right) 
\end{array}%
\right) ,
\end{equation*}%
where 
\begin{eqnarray*}
D_{n}\left( x,y\right)  &=&\frac{\partial }{\partial x}S_{n}\left(
x,y\right) , \\
\widetilde{I}_{n}\left( x,y\right)  &=&-\int_{x}^{y}S_{n}\left( x,z\right)
dz+\frac{1}{2}\mathrm{sgn}\left( y-x\right) ,
\end{eqnarray*}%
and $S_{n}\left( x,y\right) $ can be defined in terms of skew-orthogonal or
orthogonal polynomials. (See Propositions 6.3.2, 6.3.3, and 6.4.3 in \cite%
{forrester10}.)

The expected number of eigenvalues in an interval $I$ for the Jacobi
ensemble can be computed as 
\begin{equation*}
\int_{I}S_{p_{N}}\left( x,x\right) dx.
\end{equation*}%
The scaling limit for $S_{n}\left( x,y\right) $ at the edge of the spectrum
is given in Section 7.9.1 of \cite{forrester10}. We change the formula a bit
since we use interval $\left[ 0,1\right] $ instead of $[-1,1]$ 
\begin{eqnarray*}
\lim_{n\rightarrow \infty }\frac{1}{2n^{2}}S_{n}\left( \frac{x^{2}}{2n^{2}},%
\frac{y^{2}}{2n^{2}}\right) &=&:S^{hard}\left( x^{2},y^{2}\right) =\frac{x}{y%
}\left. K^{hard}\left( x,y\right) \right. _{a\rightarrow a+1} \\
&&+\frac{J_{a+1}\left( y\right) }{2y}\left( 1-\int_{0}^{x}J_{a+1}\left(
v\right) dv\right) ,
\end{eqnarray*}%
where%
\begin{equation*}
K^{hard}\left( x^{2},y^{2}\right) =\frac{xJ_{a+1}\left( x\right) J_{a}\left(
y\right) -yJ_{a+1}\left( y\right) J_{a}\left( x\right) }{\left(
x^{2}-y^{2}\right) }.
\end{equation*}%
In particular, 
\begin{equation*}
S^{hard}\left( x^{2},x^{2}\right) =\frac{1}{2}\left( J_{a+1}\left(
x^{1/2}\right) ^{2}-J_{a+2}\left( x\right) J_{a}\left( x\right) \right) +%
\frac{J_{a+1}\left( x\right) }{2x}\left( 1-\int_{0}^{x}J_{a+1}\left(
v\right) dv\right) .
\end{equation*}%
Hence, we can write 
\begin{eqnarray*}
\mathbb{E}\mathcal{N}\left( I_{N}\right)
&=&\int_{0}^{t_{0}^{2}/(2p_{N}^{2})}S_{p_{N}}\left( x,x\right) dx \\
&=&\int_{0}^{t_{0}}\frac{t}{p_{N}^{2}}S_{p_{N}}\left( \frac{t^{2}}{2p_{N}^{2}%
},\frac{t^{2}}{2p_{N}^{2}}\right) dt \\
&=&\int_{0}^{t_{0}}t\left( J_{a+1}\left( t\right) ^{2}-J_{a+2}\left(
t\right) J_{a}\left( t\right) \right) dt \\
&&+\int_{0}^{t_{0}}J_{a+1}\left( t\right) \left( 1-\int_{0}^{t}J_{a+1}\left(
v\right) dv\right) dt+o(1).
\end{eqnarray*}%
$\square $

\section{Conclusion}

\label{section_conclusion}

We showed that the local behavior of eigenvalues of the sum of two random
orthogonal projections satisfy the universality paradigm in general, with a
single exception. The exception can perhaps be explained by an additional
symmetry present in the model.

This evidence give some support to the conjecture that the local behavior of
eigenvalues of the sum of two non-scalar large random matrices is universal.

\bibliographystyle{plain}
\bibliography{comtest}

\end{document}